# APPROXIMATE ZERO-ONE LAWS AND SHARPNESS OF THE PERCOLATION TRANSITION IN A CLASS OF MODELS INCLUDING TWO-DIMENSIONAL ISING PERCOLATION[1]

By J. van den Berg

*CWI and Vrije Universiteit*

One of the most well-known classical results for site percolation on the square lattice is the equation $p_c + p_c^* = 1$. In words, this equation means that for all values $\neq p_c$ of the parameter $p$, the following holds: either a.s. there is an infinite open cluster or a.s. there is an infinite closed "star" cluster. This result is closely related to the percolation transition being sharp: below $p_c$, the size of the open cluster of a given vertex is not only (a.s.) finite, but has a distribution with an exponential tail. The analog of this result has been proven by Higuchi in 1993 for two-dimensional Ising percolation (at fixed inverse temperature $\beta < \beta_c$) with external field $h$, the parameter of the model.

Using sharp-threshold results (approximate zero-one laws) and a modification of an RSW-like result by Bollobás and Riordan, we show that these results hold for a large class of percolation models where the vertex values can be "nicely" represented (in a sense which will be defined precisely) by i.i.d. random variables. We point out that the ordinary percolation model obviously belongs to this class and we also show that the Ising model mentioned above belongs to it.

**1. Introduction.** A landmark in the development of percolation theory is the proof by Kesten in 1980 ([18]) that the critical probability for bond percolation on the square lattice equals 1/2. A key argument in his proof is what would now be called a "sharp-threshold result:" he showed that if $n$ is large and the probability of having an open horizontal crossing of an $n \times n$ box is neither close to 0, nor close to 1, then there is a reasonable probability to have many (in fact, at least of order $\log n$) so-called pivotal edges (or cut edges). [These are edges $e$ with the property that changing

Received July 2007; revised November 2007.
[1]Supported in part by the Dutch BSIK/BRICKS project.

the state (open/closed) of $e$ changes the occurrence or nonoccurrence of an open horizontal crossing.]

The proof of this intermediate key result has a combinatorial-geometric flavor: it involves a "counting argument" with conditioning on the lowest open crossing (and the leftmost closed dual crossing of the area above the open crossing just mentioned). This result, in turn, implies that the derivative (w.r.t. the parameter $p$) of the crossing probability is very large (also at least of order $\log n$) if $n$ is very large. Since probabilities are at most 1, it is impossible to have such behavior for all $p$ throughout some interval of nonzero length. On the other hand, other arguments show that the above mentioned crossing probability is bounded away from 0 and 1, uniformly in $n$ and $p \in (1/2, p_c)$. Hence $p_c$ must be equal to $1/2$. [It was already known ([14]) that $p_c \geq 1/2$.]

Soon after Kesten's result, it was shown ([28] and [32]) that his arguments can also be used to prove related long-standing conjectures, in particular, that $p_c + p_c^* = 1$ for site percolation on the square lattice. Here, $p_c^*$ denotes the critical probability for the so-called matching lattice (or star lattice) of the square lattice—this is the lattice with the same vertices as the square lattice, but where each vertex $(x, y)$ has not only edges to its four horizontal or vertical "neighbors" $\{(x', y') : |x - x'| + |y - y'| = 1\}$, but also to the nearest vertices in the diagonal directions, $\{(x', y') : |x - x'| = |y - y'| = 1\}$.

Russo [29] was the first to put the aforementioned "sharp-threshold" argument of Kesten in a more general framework by formulating an *approximate zero-one law*. This approximate zero-one law is *not* itself a percolation result. It is, as the name indicates, a "finite" approximation of Kolmogorov's zero-one law. Recall that the latter says (somewhat informally) that if $X_1, X_2, \ldots$ are i.i.d. Bernoulli random variables (say, with parameter $p$) and $A$ is an event with the property that the occurrence or nonoccurrence of $A$ cannot be changed by changing a single $X_i$, then $A$ has probability 0 or 1. Russo's approximate law says that if $A$ is an increasing event with the property that for each $i$ and $p$, the probability that changing the state of $X_i$ disturbs the occurrence or nonoccurrence of $A$ is very small, then for all $p$, except on a very small interval, the probability of $A$ is close to 0 or close to 1.

Given this approximate zero-one law, the combinatorial-geometric argument in Kesten's work discussed above [to get a lower bound for the (expected) number of pivotal items] can be (and was, in Russo's paper) replaced by the considerably simpler (and "smoother") argument that the probability that a given edge (or, for site percolation, site) is pivotal is small when $n$ is large.

It should be noted that for the other, more standard, part of the proof, the approximate zero-one law does not help: a "separate" argument of the form that if the probability of an open crossing of a square is sufficiently large, then there is an infinite open cluster is needed. For this argument, which is often called a "finite-size criterion," the RSW theorem ([27] and [30]; see also [13], Chapter 11) is essential. Informally, this theorem gives a suitable lower bound for the probability of having an open crossing (in the "long" direction) of a $2n \times n$ rectangle, in terms of the probability of having a crossing of an $n \times n$ square. The classical proof uses conditioning on the lowest crossing. For ordinary Bernoulli percolation, this works fine, but, as remarked earlier, in dependent models, such conditioning often leads to very serious, if not unsolvable, problems. An important recent achievement in this respect is a "box-crossing" theorem obtained in [8], the proof of which is much more robust than that of the "classical" RSW result and does not use such conditioning. We will come back to this later.

Refinements, generalizations and independent results with partly the same flavor as Russo's approximate zero-one law have been obtained and/or applied by Kahn, Kalai and Linial [17], Talagrand [31], Friedgut and Kalai [10], Bollobás and Riordan (see, e.g., [7]), Graham and Grimmett [12], Rossignol [25] and others, and have become known as "sharp-threshold results."

Although such results take Kesten's key argument in a more general context, involving "less geometry," Kesten's proof is still essentially the shortest and, from a probabilistic point of view, intuitively the most appealing, self-contained proof of $p_c = 1/2$ for bond percolation on the square lattice: none of the aforementioned general sharp-threshold theorems has a short or probabilistically intuitive proof. Moreover, the combinatorial-geometric ideas and techniques in Kesten's proof have turned out to be very fruitful in other situations, for instance, in the proof of one of the main results in Kesten's paper on scaling relations for two-dimensional percolation ([20]).

On the other hand, there are examples of percolation models where Kesten's argument is difficult and cumbersome to carry out, or where it is even not (yet) known how to do this. An example of the latter is the Voronoi percolation model, for which Bollobás and Riordan ([8]; see also [9]) proved (using a sharp-threshold result from [10]) that it has critical probability $1/2$. This had been conjectured for a long time, but carrying out Kesten's strategy for that model led to (thus far) unsolved problems.

An example of the former is percolation of $+$ spins in the two-dimensional Ising model with fixed inverse temperature $\beta < \beta_c$ and external field parameter $h$ (which plays the role of $p$ in "ordinary" percolation). Higuchi ([15] and [16]) showed that for all values of $h$, except the critical value $h_c$, either (a.s.) there is an infinite cluster of vertices with spin $+$ or (a.s) there is an infinite $*$ cluster (i.e., a cluster in the $*$ lattice) of vertices with spin $-$. (The result is stated in [16], but much of the work needed in the proof is done in [15].) Higuchi followed the scheme of Kesten's arguments. However, to carry them out (in particular the "counting under conditioning on the lowest crossing," etc.), he had to overcome several new technical difficulties due to the dependencies in this model. This makes the proof far from easy to read.

REMARK. It should be mentioned here that there is also a very different proof of $p_c = 1/2$ for bond percolation on $\mathbb{Z}^2$ (and $p_c + p_c^* = 1$ for site percolation), namely, by using the work of Menshikov [23] and of Aizenman and Barsky [1]. They gave a more "direct" (not meaning "short" or "simple") proof, without using the results or arguments indicated above, that for every $d \geq 1$, the cluster radius distribution for independent percolation on $\mathbb{Z}^d$ with $p < p_c$ has an exponential tail. However, their proofs use (and need) the BK inequality ([5]) and since our interest is mainly in dependent percolation models (for which no suitable analog of this inequality seems to be available), these proof methods will not be discussed in more detail here.

In the current paper, we present a theorem (Theorem 2.2) which says that the analog of $p_c + p_c^* = 1$ holds for a large class of weakly dependent two-dimensional percolation models. Roughly speaking, this class consists of systems that have a proper, monotone, "finitary" representation in terms of i.i.d. random variables. From the precise definitions, it will be immediately clear that it contains the ordinary (Bernoulli) percolation models. We give (using results obtained in the early 1990s by Martinelli and Olivieri [22] and modifications of results in [6] which were partly inspired by [24]) a "construction" of the earlier mentioned two-dimensional Ising model which shows that this model also belongs to this class.

We hope that our results improve insight into the Ising percolation model and will help to show that many other (not yet analyzed) weakly dependent percolation models also belong to the aforementioned class.

The theorem is based on:

(a) One of the sharp-threshold results mentioned above, namely Corollaries 1.2 and 1.3 in [31] (and a recent generalization in [26]), which are close in spirit to, but quantitatively more explicit than, Russo's approximate zero-one law. The reason for using these sharp-threshold results rather than those in [10] (which, as stated above, were applied to percolation problems by Bollobás and Riordan) is that the latter assume certain symmetry properties on the events to which they are applied. In many situations, this causes no essential difficulties, but it gives much more flexibility to allow an absence of such symmetries (see the remark following property (iv) near the end of Section 2.1).

(b) A modification/improvement (obtained in [4]) of an RSW-like box-crossing theorem of Bollobás and Riordan [8]. As indicated in the short discussion of Russo's paper above, some form of RSW theorem seems unavoidable. For many dependent percolation models, it is very hard (or maybe impossible) to carry out the original proof of RSW. The Bollobás–Riordan form of RSW (and its modification in [4]) is weaker (but still strong enough) and much more robust with respect to spatial dependencies.

In Section 2, we introduce some terminology and state Theorem 2.2, which says that a large class of two-dimensional percolation models satisfies an analog of $p_c + p_c^* = 1$. We also state some consequences/examples of the theorem. In particular, we show that the Ising percolation model studied by Higuchi satisfies the conditions of Theorem 2.2 so that his result mentioned above can be alternatively derived from our theorem.

In Section 3, we state preliminaries needed in the proof of Theorem 2.2: Talagrand's result mentioned above and an extension of his result to the case where the underlying random variables can take more than two different values, and where the events under consideration do not necessarily depend on only finitely many of these underlying variables. In that section, we also explain that the earlier mentioned (modification of the) RSW-result of Bollobás and Riordan applies to our class of percolation models and we prove other properties that are used in the proof of Theorem 2.2.

In Section 4, we finally prove Theorem 2.2, using the ingredients explained in Section 3.

Apart from the proofs of the RSW-like theorem and of Talagrand's sharp-threshold result mentioned above, the proof of Theorem 2.2 is practically self-contained.

## 2. Statement of the main theorem and some corollaries.

2.1. *Terminology and set-up.* In this subsection, we will describe the (dependent) percolation models on the square lattice for which our main result, a generalization of the well-known $p_c + p_c^* = 1$ for ordinary percolation, holds.

Throughout this paper, we use the norm

$$\|v\| := |v_1| + |v_2|,$$

where $v = (v_1, v_2) \in \mathbb{Z}^2$.

Let $k$ be a positive integer and let $\mu^{(h)}$, $h \in \mathbb{R}$, be a family of probability measures on $\{0, 1, \ldots, k\}$, indexed by the parameter $h$, with the following two properties:

(a) for each $1 \leq j \leq k$, $\mu^{(h)}(\{j, \ldots, k\})$ is a continuously differentiable, strictly increasing function of $h$;

(b) $\lim_{h \to \infty} \mu^{(h)}(k) = \lim_{h \to -\infty} \mu^{(h)}(0) = 1$.

Let $I$ be a countable set. Before we go on, we need some notation and a definition. We use "$\subset\subset$" to indicate "finite subset of." The special elements $(0, 0, 0, \ldots)$ and $(k, k, k, \ldots)$ of $\{0, \ldots, k\}^I$ are denoted by **0** and **k**, respectively.

Let $f:\{0,\ldots,k\}^I \to \mathbb{R}$ be a function. Let $V \subset\subset I$ and let $y \in \{0,\ldots,k\}^V$. For $x \in \{0,\ldots,k\}^I$, we write $x_V$ for the "tuple" $(x_i, i \in V)$. We say that $y$ determines (the value of) $f$ if $f(x) = f(x')$ for all $x, x'$ with $x_V = x'_V = y$.

Let $X_i$, $i \in I$, be independent random variables, each with distribution $\mu^{(h)}$. Let $P^{(h)}$ denote the joint distribution of the $X_i$'s. The $X_i$'s will be the "underlying" i.i.d. random variables for our percolation system. We will assume that the "actual spin variables," which take values $+1$ ("open") and $-1$ ("closed") and which will be denoted by $\sigma_v, v \in \mathbb{Z}^2$ below, are "suitably described" in terms of the underlying $X$ variables: for each $v \in \mathbb{Z}^2$, its spin variable $\sigma_v$ is a function of the $(X_i, i \in I)$. These functions do not themselves depend on $h$, but changing $h$ will change the distribution of the $X$ variables and thus that of the $\sigma$ variables. More precisely, we assume that $\sigma_v, v \in \mathbb{Z}^2$, are random variables with the following properties:

(i) (*Monotonicity.*) For each $v$, $\sigma_v$ is a measurable, increasing, $\{-1,+1\}$-valued function of the collection $(X_i, i \in I)$ and, moreover, for each $v \in \mathbb{Z}^2$, $\sigma_v(\mathbf{0}) = -1$ and $\sigma_v(\mathbf{k}) = +1$;

(ii) (*Finitary representation.*) There exist $C_0 > 0$ and $\gamma > 0$ such that for each $v \in \mathbb{Z}^2$, there is a sequence $i_1(v), i_2(v), \ldots$ of elements of $I$ such that for all positive integers $m$ and all $h \in \mathbb{R}$,

$$P^{(h)}((X_{i_1(v)},\ldots,X_{i_m(v)}) \text{ does not determine } \sigma_v) \leq \frac{C_0}{m^{2+\gamma}};$$

(iii) (*Mixing*).

(1)
$$\exists \alpha > 0 \ \forall v,w \in \mathbb{Z}^2 \ \forall m < \alpha \|v-w\|,$$
$$\{i_1(v),\ldots,i_m(v)\} \cap \{i_1(w),\ldots,i_m(w)\} = \varnothing;$$

(iv) for each $h$, the distribution of $(\sigma_v, v \in \mathbb{Z}^2)$ is translation invariant and invariant under rotations by 90 degrees, and under vertical and horizontal axis reflection.

REMARK. Note that in property (iv), we do not require that we can identify $I$ with $\mathbb{Z}^2$ in such a way that there is a stationary mapping (i.e., a mapping which commutes with shifts) from the process $(X_i, i \in \mathbb{Z}^2)$ to the process $(\sigma_v, v \in \mathbb{Z}^2)$. In many cases, there will be such identification, but we found its requirement unnecessarily strong for our purposes (see, e.g., the example of the Ising model below, where the mapping under consideration is not of this form). A consequence of the absence of this requirement is an absence of certain symmetries needed to apply the sharp-threshold results in, for example, [10]. This is the main reason for using the results in [31].

DEFINITION 2.1. If a random field $(\sigma_v, v \in \mathbb{Z}^2)$ has the properties (i)–(iv) above, we say that the process has a *nice, finitary representation* (in terms of the $X$ process and with parameter $h \in \mathbb{R}$).

2.2. *Statement of the main theorem and some special cases.* Now we consider percolation in terms of the $\sigma$ variables: we interpret $\sigma_v = +1$ (resp. $-1$) as the vertex $v$ being open (resp. closed) and are interested in (among other things) the existence of infinite paths on which every vertex is open. As usual, in our notion of "ordinary" paths, we allow only horizontal and vertical steps and we use the term "star paths" when, in addition to these steps diagonal steps are also allowed. Similarly (and, again, following the usual conventions), we define "ordinary" clusters as well as star clusters. When we speak simply of a "cluster," we will always mean an "ordinary" cluster.

The + cluster of a vertex $v$ will be denoted by $C_v^+$; the $-*$ cluster of $v$ (i.e., the $-$ cluster of $v$ in the star lattice) will be denoted by $C_v^{-*}$, etc. If $v = 0$ [the vertex $(0,0)$], we will often omit the subscript $v$.

Recall that $P^{(h)}$ denotes the probability distribution of the collection $(X_i, i \in I)$. We will also use it for the probability measure on $\{-1, +1\}^{\mathbb{Z}^2}$ induced by the map from the $X$ variables to the $\sigma$ variables. Since the context in which it is used will always be clear, this should not cause any confusion.

THEOREM 2.2.  *Let $(\sigma_v, v \in \mathbb{Z}^2)$ be a spin system with a nice, finitary representation, with parameter $h \in \mathbb{R}$ (in the sense of Definition 2.1).*

*Then there is a critical value $h_c$ of $h$ such that:*

(a) $\forall h > h_c P^{(h)}(|C^+| = \infty) > 0$ *and the distribution of $|C^{-*}|$ has an exponential tail;*

(b) $\forall h < h_c P^{(h)}(|C^{-*}| = \infty) > 0$ *and the distribution of $|C^+|$ has an exponential tail.*

REMARK.  Note that it follows from the statement of the theorem that $h_c$ satisfies
$$h_c = \inf\{h : P^{(h)}(|C^+| = \infty) > 0\} = \sup\{h : P^{(h)}(|C^{-*}| = \infty) > 0\}.$$

Also, note that if reversal of $h$ corresponds with a spin-flip [more precisely, if the distribution of $\sigma$ under $P^{(h)}$ is the same as the distribution of $-\sigma$ $(= (-\sigma_v, v \in \mathbb{Z}^2))$ under $P^{(-h)}$], the above theorem immediately implies that

(2) $\quad h_c + h_c^* = 0, \qquad$ where $h_c^* = \inf\{h : P^{(h)}(|C^{+*}| = \infty) > 0\}.$

2.2.1. *Special cases.*

*Bernoulli site percolation on the square lattice, with parameter $p$.* This model, where the vertices are open $(+1)$ with probability $p$ and closed $(-1)$ with probability $1-p$ trivially satisfies the conditions of Theorem 2.2: simply take $I = \mathbb{Z}^2$, $k = 1$ (i.e., the $X_i$'s take values 0 and 1) and $\sigma_v = 2X_v - 1$, $v \in \mathbb{Z}^2$. Finally, take, for instance, [note that we want $\mu^{(h)}(1)$ to go 1 (resp. 0) as $h \to \infty$ $(-\infty)$]

$$\mu^{(h)}(1) = \frac{\exp(h)}{\exp(h) + \exp(-h)}.$$

Taking $p = \mu^{(h)}(1)$ completes the "translation." It is easy to see that reversing $h$ corresponds with a spin-flip, so (2) holds, which is equivalent to the well-known

$$p_c + p_c^* = 1$$

for this model.

*Models defined explicitly in terms of i.i.d. random variables.* In the previous example, the representation in terms of i.i.d. random variables was explicit and trivial. It is easy to find many other examples with explicit (but less trivial) representations. For instance, take $I = \mathbb{Z}^2$ and let the $X$ variables be i.i.d. Bernoulli with parameter $p$. Define, for each $v \in \mathbb{Z}^2$, $\sigma_v$ as follows. Consider, for each $n$, the difference between the number of 1's and the number of 0's in the $2n \times 2n$ square centered at $v$. Take the smallest $n$ where this difference has absolute value larger than some constant, say 5. Define $\sigma_v$ as the sign of the aforementioned difference (number of 1's minus number of 0's) for that $n$. It is easy to check that this definition corresponds to a nice, finitary representation, in the sense of Definition 2.1. More interesting (in the context of the subject of this paper) are those weakly dependent models that are not a priori explicitly defined in terms of such a representation. One can then search for a possible "hidden" representation. A major example where this works is the following.

*Ising model with (fixed) inverse temperature $\beta < \beta_c$ and external field parameter $h$.* We first recall some definitions and standard results for these models. Ising measures $\mu_{\beta,h}$ on $\{-1,+1\}^{\mathbb{Z}^2}$, with inverse temperature $\beta \in [0,\infty)$ and external field $h \in (-\infty,\infty)$, are probability measures that satisfy, for $\eta \in \{0,1\}$ and $v \in \mathbb{Z}^2$,

$$
\begin{aligned}
&\mu_{\beta,h}(\sigma_v = \eta \mid \sigma_w, w \neq v) \\
&\quad = \frac{\exp(\beta\eta(h + \sum_{w \sim v} \sigma_w))}{\exp(\beta\eta(h + \sum_{w \sim v} \sigma_w)) + \exp(-\beta\eta(h + \sum_{w \sim v} \sigma_w))},
\end{aligned}
\tag{3}
$$

where $w \sim v$ means that $\|v - w\| = 1$.

It is well known that there is a critical value $\beta_c$ such that for $\beta < \beta_c$, there is a unique measure satisfying (3), while for $\beta > \beta_c$, there is more than one such measure.

The Ising model is one of the most well-known examples of a Markov random field: the conditional distribution of the spin value of a vertex $v$, given the spin values of all other vertices, depends only on the spin values of the neighbors of $v$.

The "single-site" conditional distributions in (3) will often be used in the remainder of this subsection and will be denoted by $q_v^\alpha$. More precisely, for $v \in \mathbb{Z}^2$, let $\partial v$ denote the set of (four) vertices that are neighbors of $v$. Further, for $\alpha \in \{-1,+1\}^{\mathbb{Z}^2}$ and $V \subset\subset \mathbb{Z}^2$, let $\alpha_V$ denote the "restriction" of $\alpha$ to $V$; that is, $\alpha_V = (\alpha_w, w \in V)$. For $\alpha \in \{-1,+1\}^{\partial v}$ and $\eta \in \{-1,+1\}$, we define $q_v^\alpha(\eta) = q_v^\alpha(\eta; \beta, h)$ as the conditional probability that $\sigma_v$ equals $\eta$, given that $\sigma_{\partial v}$ equals $\alpha$:

$$(4) \qquad q_v^\alpha(\eta) := \frac{\exp(\beta\eta(h + \sum_{w \sim v} \alpha_w))}{\exp(\beta\eta(h + \sum_{w \sim v} \alpha_w)) + \exp(-\beta\eta(h + \sum_{w \sim v} \alpha_w))}.$$

Note that the dependence on the "neighbor configuration" $\alpha$ is only through the *number* of $+$ (and hence of $-$) spins in $\alpha$. Therefore, it is also convenient to define, for $m = 0, \ldots, 4$,

$$(5) \qquad q_v^{(m)}(\eta) = q_v^\alpha(\eta),$$

where $\alpha$ may be any element of $\{-1,+1\}^{\partial v}$ with the property that the number of $w \sim v$ with $\alpha_w = +1$ equals $m$.

The following result is well known and goes back to [2] and [21]. Higuchi [15] proved and used a stronger result, but the weaker version below is sufficient for our purposes.

THEOREM 2.3. *There exist $C_1 > 0$ and $\lambda_1 > 0$ (which depend on $\beta$, but not on $h$) such that*

$$(6) \qquad \begin{aligned} &\mu_{\beta,h}(\sigma_0 = +1 \mid \sigma_{\partial\Lambda(n)} \equiv +1) - \mu_{\beta,h}(\sigma_0 = +1 \mid \sigma_{\partial\Lambda(n)} \equiv -1) \\ &\leq C_1 \exp(-\lambda_1 n). \end{aligned}$$

*Here, $\Lambda(n)$ denotes the set of vertices $[-n,n]^2$ and $\partial\Lambda(n)$ the boundary of this set.*

Martinelli and Olivieri (Theorem 3.1 in [22]) have proven, for a large class of spin systems, that such a *spatial* mixing property implies exponential convergence (to equilibrium) for certain dynamics. For the Ising model, this dynamics is as follows. First, we define the notion *local update*. Let $\alpha \in \{-1,+1\}^{\mathbb{Z}^2}$ and $v \in \mathbb{Z}^2$. By a local update of the spin value of $v$ (in the configuration $\alpha$), we mean that we draw a new value, say $\eta$, according to the distribution $q_v^{\alpha_{\partial v}}(\cdot)$ and leave $\alpha$ unchanged everywhere, except at $v$, where we replace $\alpha_v$ by $\eta$. The dynamics can now be described as follows. Start from some initial configuration. Each vertex is *activated* at rate 1. When a vertex is activated, a local update at that vertex is made. The Martinelli–Olivieri result (for the special case of the Ising model) says that the distribution at time $t$, starting from any initial configuration, converges exponentially fast (uniformly in $h$) to $\mu_{\beta,h}$. In particular, the probability that 0 has spin value $+1$ at time $t$ converges exponentially fast (and uniformly in $h$ and in the initial configuration) to $\mu_{\beta,h}(\sigma_0 = +1)$.

As observed in [6], this also holds for certain discrete-time versions of the dynamics. The discrete-time dynamics in [6] involves auxiliary random variables, in terms of which the dynamics is not monotone. For the purposes of [6], that did not matter, but this dynamics is not suitable for our current purpose—to "construct" the Ising measure in such a way that it fits with Definition 2.1. The following dynamics *is* suitable for our purposes and the Martinelli–Olivieri proof (with straightforward modifications) works for this dynamics as well. In this discrete-time dynamics, we update all even vertices at the even times and all odd vertices at the odd times. [A vertex is even (resp. odd) if the sum of its coordinates is even (resp. odd).] Note that these "parallel" updates are well defined since the update of an even (resp. odd) vertex only involves the "current" spin values of its neighbors, each of which is odd (resp. even).

To describe the Ising model as a nice, finitary representation, in the sense of Definition 2.1, we describe these local updates as follows in terms of i.i.d. random variables $Y_i(t), i \in \mathbb{Z}^2, t \in \mathbb{N}$, which take values in $\{-1, 0, \ldots, 4\}$. Here (and further), $\sigma_v^\omega(t)$ denotes the spin value at vertex $v$ at time $t$ for the system starting at time 0 with configuration $\omega$. Sometimes, we will omit the superscript $\omega$. At each even time $t$, we do the following, for each even vertex $v$: if the number of $w \sim v$ with $\sigma_t(w) = -1$ is at most $Y_v(t)$, we set $\sigma_v(t+1) := +1$, otherwise we set $\sigma_v(t+1) := -1$. For odd $t$, we perform the analogous actions for all odd $v$. It is easy to see [recall (5)] that if we take the following distribution for the $Y$ variables, these actions correspond exactly with the earlier defined notion of local updates:

$$P(Y_v(t) \geq m) = q_0^{(4-m)}(+1; \beta, h), \qquad 0 \leq m \leq 4.$$

An advantage of using such auxiliary variables is that it enables the coupling of systems starting from different initial configurations. Define $\sigma^\omega(t) = (\sigma_v^\omega(t), v \in \mathbb{Z}^2)$ as the configuration at time $t$ for the system that starts at time 0 with configuration $\omega$ and follows the aforementioned dynamics (involving the $Y$ variables). We will simply replace the superscript $\omega$ by $+$ when we start with the initial configuration where each vertex has value $+1$, and by $-$ when we start with $-$ values. As said before, the Martinelli–Olivieri result concerning exponential convergence to equilibrium extends to this dynamics. In terms of the above notation, this Martinelli–Olivieri result tells us that there are positive $C_2$ and $\lambda_2$ (which depend on $\beta$ but not on $h$) such that for all $t$,

(7) $$P(\sigma_v^+(t) \neq \sigma_v^-(t)) \leq C_2 \exp(-\lambda_2 t).$$

Also, note that we can extend the collection of $Y$ variables to negative $t$ and that for all integers $s, t$ with $s \leq t$, and all configurations $\omega \in \{-1, +1\}^{\mathbb{Z}^2}$, we can define $\sigma^\omega(s,t) = (\sigma_v^\omega(s,t), v \in \mathbb{Z}^2)$ as the configuration at time $t$ for the system that starts at time $s$ with configuration $\omega$ and evolves as described above. Analogously as in [6] (which was partly inspired by the perfect simulation ideas in [24]), we observe that if $t < 0$ and $\sigma_v^+(t,0) = \sigma_v^-(t,0)$, then (by obvious monotonicity) $\sigma_v^\omega(s,0) = \sigma_v^{\omega'}(s,0)$ for all $s \leq t$ and all $\omega, \omega'$. From this observation, (7) and standard arguments, it follows that if we define

$$\tau(v) = \max\{t < 0 : \sigma_v^+(t,0) = \sigma_v^-(t,0)\}, \qquad v \in \mathbb{Z}^2$$

and

(8) $$\sigma(v) = \sigma_v^+(\tau(v), 0) \qquad (= \sigma_v^-(\tau(v), 0)), \qquad v \in \mathbb{Z}^2,$$

then we have that $\sigma := (\sigma(v), v \in \mathbb{Z}^2)$ has the Ising distribution $\mu_{\beta, h}$ and that

(9) $$P(\tau(v) \geq n) \leq C_2 \exp(-\lambda_2 n).$$

This shows that the Ising distribution indeed has a nice, finitary representation (in the sense of Definition 2.1). Take $I = \{(v,t) : v \in \mathbb{Z}^2, t \in \mathbb{Z}, t < 0\}$ and $X_{(v,t)} = Y_v(t), (v,t) \in I$. Then (i) is clear. To see (ii), note that for each $t < 0$, $\sigma_v^+(t,0)$ and $\sigma_v^-(t,0)$ are completely determined by the variables $Y_w(s)$, $t \leq s < 0$, $\|w - v\| < s$.

So, for the sequence $i_1(v), i_2(v), \ldots$, we can take $(v, -1)$, followed by an enumeration of the (finite) set $\{(w, -2) : w \in \mathbb{Z}^2, \|w - v\| < 2\}$, followed by an enumeration of $\{(w, -3) : w \in \mathbb{Z}^2, \|w - v\| < 3\}$, etc. The upper bound in (ii) (in fact, even a stronger bound) for the probability that $X_{i_1(v)}, \ldots, X_{i_m(v)}$ does not determine $\sigma_v$ follows from (9) and the fact that the set $\{(w,s) : \|w - v\| < |s|, t \leq s < 0\}$ has of order $|t|^3$ elements. Property (iii) is now also clear. Property (iv) is standard (and has nothing to do with the above description of the Ising model in terms of the $Y$ variables: since $\beta < \beta_c$, there is a unique Ising measure with parameters $\beta, h$ and this measure inherits the symmetry properties in the definition of the model).

Hence, we may apply Theorem 2.2. Moreover, the spin-flip symmetry mentioned in the remark preceding (2) is clearly satisfied. So, we get the following, which is the result by Higuchi mentioned earlier (see Theorem 1 (and Corollary 2) in [16]).

THEOREM 2.4. *Let $\beta < \beta_c$ and consider the Ising measures $\mu_{\beta, h}$, $h \in \mathbb{R}$, on the square lattice. Statements* (a) *and* (b) *of Theorem 2.2 above (with $P^{(h)} = \mu_{\beta, h}$), as well as equation (2), hold for this model.*

REMARKS. (i) The sharp-threshold result in [12] may provide yet another route to prove this result for the Ising model. However, that sharp-threshold result is not suitable for the proof of our general Theorem 2.2 because the random field $\sigma_i, i \in \mathbb{Z}^2$ in Theorem 2.2 does not necessarily satisfy the strong FKG condition needed in [12].

(ii) We hope that, like the Ising model, many other models which at first sight are not covered by Theorem 2.2 can be constructed or represented in such a way that this theorem does apply. However, we do not claim that this theorem gives a completely general recipe. For instance, attempts to bring the models treated in [3] (which have some of the flavor of the Ising model) into the context of this theorem have, thus far, not been successful.

## 3. Preliminaries.

3.1. *Approximate zero-one laws.* A key ingredient in our proof of Theorem 2.2 is a sharp-threshold result (or approximate zero-one law). As stated in Section 1, there are several such results in the literature. The one we use is Corollary 1.2 from Talagrand's paper [31], which is somewhat similar in spirit to Russo's approximate zero-one law ([29]), but more (quantitatively) explicit.

These threshold results are, although particularly useful for percolation, of a much more general nature. Consider the set $\Omega := \{0,1\}^n$, which, for the moment, serves as our sample space. For $\omega, \omega' \in \Omega$, we say that $\omega \leq \omega'$ (or, equivalently, $\omega' \geq \omega$) if $\omega_i \leq \omega_i'$ for all $1 \leq i \leq n$. Following the standard terminology, we say that an event (subset of $\Omega$) is increasing if for each $\omega \in A$ and each $\omega' \geq \omega$, we have $\omega' \in A$. For $\omega \in \Omega$ and $1 \leq i \leq n$, we define $\omega^{(i)}$ as the configuration obtained from $\omega$ by flipping $\omega_i$. More precisely, $\omega_j^{(i)}$ is equal to $\omega_j$ for $j \neq i$, and $1 - \omega_j$ if $j = i$.

Let $A$ be an increasing event, $\omega \in \Omega$ and $1 \leq i \leq n$. We say that $i$ is an internal pivotal index (for $A$, in the configuration $\omega$) if $\omega \in A$, but $\omega^{(i)} \notin A$. It is easy to see from the fact that $A$ is increasing that this implies that $\omega_i = 1$.

By $A_i$, we denote the event that $i$ is an internal pivotal for $A$; that is,

$$A_i = \{\omega : \omega \in A \text{ but } \omega^{(i)} \notin A\}.$$

Let, for $p \in (0,1)$, $\mathbb{P}_p$ be the product measure with parameter $p$. Talagrand's result to which we referred above is the following.

THEOREM 3.1 ([31], Corollary 1.2). *There is a universal constant $K_1$ such that for all $n$, all increasing events $A \subset \{0,1\}^n$ and all $p$,*

$$\frac{d}{dp} \mathbb{P}_p(A) \geq \frac{\log(1/\varepsilon)}{K_1} \mathbb{P}_p(A)(1 - \mathbb{P}_p(A)), \tag{10}$$

*where $\varepsilon = \varepsilon(p) = \sup_{i \leq n} \mathbb{P}_p(A_i)$.*

REMARK. In fact, Corollary 1.2 in [31] is somewhat sharper, namely with $K_1$ above replaced by $Kp(1-p)\log[(2/(p(1-p))]$, where $K$ is also a universal constant. Since $p(1-p)\log[(2/(p(1-p))]$ is bounded from above, Theorem 3.1 above follows immediately.

Let $p_1 < p_2$. Noting (as in Section 3 of [31]) that (10) is equivalent to

$$\frac{d}{dp}\log\left(\frac{\mathbb{P}_p(A)}{1-\mathbb{P}_p(A)}\right) \geq \frac{\log(1/\varepsilon)}{K_1}$$

and integrating this inequality over the interval $(p_1, p_2)$ gives the following.

COROLLARY 3.2 ([31], Corollary 1.3). *There is a universal constant $K_1$ such that for all $n$, all increasing events $A \subset \{0,1\}^n$ and all $p_1 < p_2$,*

(11) $$\mathbb{P}_{p_1}(A)(1 - \mathbb{P}_{p_2}(A)) \leq (\varepsilon')^{(p_2-p_1)/K_1},$$

*where*

(12) $$\varepsilon' = \sup_{p_1 \leq p \leq p_2} \max_{1 \leq i \leq n} \mathbb{P}_p(A_i).$$

REMARK. In the definition of $\varepsilon'$ in the statement of Corollary 1.3 in [31], the supremum involving $p$ is over the interval $[0,1]$ instead of $[p_1, p_2]$, but it is clear that the result with $\varepsilon'$ defined as in (12) holds.

We want to apply similar results to the family of distributions $P^{(h)}, h \in \mathbb{R}$, in the statement of Theorem 2.2. Recall that $P^{(h)}$ is the product over $I$ of the distribution of $\mu^{(h)}$ and that the latter is a probability distribution on $\{0, \ldots, k\}$. First, we must "generalize" some of our definitions.

The notion of increasing event is extended in the obvious way. The extension of the notion of being *pivotal* is somewhat less obvious. Let $A \subset \{0, 1, \ldots, k\}^I$ be an increasing event. We say that index $i \in I$ is an internal pivotal index (in a configuration $\omega \in \{0, \ldots, k\}^I$ and for a given increasing event $A$) if $\omega \in A$, but $\omega^{(i)} \notin A$, where, now, $\omega^{(i)}$ is defined as the configuration $\omega'$ which has $\omega'_j = \omega_j$ for all $j \neq i$ and $\omega'_i = 0$. (It follows immediately from the definition that if $i$ is pivotal, then $\omega_i > 0$.)

We cannot immediately use Corollary 3.2 because of the following two issues: one is that $k$ may be larger than 1, the other is that $I$ is not finite, but countably infinite. As to the first issue, an extension of Corollary 3.2 to $k > 1$ can be obtained by suitably "coding" $\{0, 1, \ldots, k\}$-valued random variables in terms of 0–1 valued random variables. As to the second issue, that can be overcome by restricting to a suitable subclass of increasing events (which turns out to be sufficient, but is not very satisfactory). The strategy followed by Rossignol (see [26]) is considerably more powerful. Roughly speaking, he extends Theorem 1.5 in [31] (which is a "functional" generalization of Theorem 1.1 in [31], of which Theorem 3.1 above is an easy consequence) and, from that extension, obtains the following extension of Corollary 3.2.

THEOREM 3.3 ([26], Corollary 3.1). *If the event $A \subset \{0, 1, \ldots, k\}^I$ is increasing, then for all $-\infty < h_1 < h_2 < \infty$,*

$$(13) \qquad P^{(h_1)}(A)(1 - P^{(h_2)}(A)) \leq (\bar{\varepsilon})^{(h_2-h_1)c(h_1,h_2)/K_2},$$

*where $K_2$ is a constant, $\bar{\varepsilon} = \sup_i \sup_{h \in (h_1,h_2)} P^{(h)}(A_i)$ and*

$$c(h_1, h_2) = \inf_{h \in [h_1, h_2]} \min_{1 \leq j \leq k} \frac{d}{dh} \mu^{(h)}(\{j, \ldots, k\}).$$

REMARK. In fact, Corollary 3.1 in [26] is somewhat sharper (see the remark at the end of Section 3 in [26]), but Theorem 3.3 is sufficient for our purposes.

3.2. *Mixing property.* In this subsection, we show that random variables $\sigma_v, v \in \mathbb{Z}^2$ that satisfy properties (i)–(iv) in Section 2.1 have certain, very convenient, spatial mixing properties.

We say that a vertex $v$ is $l$ *determined* (w.r.t. the $X$ configuration) if $X_{i_1}(v), \ldots, X_{i_l}(v)$ determine $\sigma_v$. A set of vertices $W$ is said to be $l$ determined if every $v \in W$ is $l$ determined. From property (ii) in Section 2.1, we have

$$(14) \qquad \begin{aligned} P^{(h)}(W \text{ not } l \text{ determined}) &\leq |W| \max_{v \in W} P^{(h)}(v \text{ not } l \text{ determined}) \\ &\leq |W| \frac{C_0}{l^{2+\gamma}} \end{aligned}$$

with $C_0$ as in property (ii).

The distance between two subsets $U, V \subset \mathbb{Z}^2$ is defined as $\min_{u \in U, v \in V} \|u - v\|$.

LEMMA 3.4. *Let $k$ be a positive integer and let $U$ and $V$ be finite subsets of $\mathbb{Z}^2$ that have distance larger than $k$ to each other. Let $A$ be an event that is defined in terms of the random variables $\sigma_v, v \in U$ and $B$ an event that is defined in terms of the random variables $\sigma_v, v \in V$. Then, with $\alpha$ and $\gamma$ as in properties* (ii) *and* (iii) *from Section 2.1,*

$$(15) \qquad |P^{(h)}(A \cap B) - P^{(h)}(A)P^{(h)}(B)| \leq 2(|U| + |V|)\frac{C_0}{\lfloor \alpha k \rfloor^{2+\gamma}}.$$

PROOF. Let $\hat{A}$ be the event $A \cap \{U \text{ is } \lfloor \alpha k \rfloor \text{ determined}\}$ and $\hat{B}$ the event $B \cap \{V \text{ is } \lfloor \alpha k \rfloor \text{ determined}\}$. Note that for each vertex $v$ and each integer $l$, the event that $v$ is $l$ determined depends only on the random variables $X_{i_1(v)}, \ldots, X_{i_l(v)}$. This, property (iii) in Section 2.1 and the fact that $U$ and $V$ have distance larger than $k$ collectively imply that $\hat{A}$ and $\hat{B}$ are independent:

$$P^{(h)}(\hat{A} \cap \hat{B}) = P^{(h)}(\hat{A})P^{(h)}(\hat{B}). \tag{16}$$

Further, using (14),

$$P^{(h)}(A \setminus \hat{A}) \leq P^{(h)}(U \text{ not } \lfloor \alpha k \rfloor \text{ determined}) \leq |U| \frac{C_0}{\lfloor \alpha k \rfloor^{2+\gamma}} \tag{17}$$

and, similarly,

$$P^{(h)}(B \setminus \hat{B}) \leq P^{(h)}(V \text{ not } \lfloor \alpha k \rfloor \text{ determined}) \leq |V| \frac{C_0}{\lfloor \alpha k \rfloor^{2+\gamma}}. \tag{18}$$

From (16)–(18), Lemma 3.4 follows straightforwardly. □

3.3. *Positive association.* The next lemma concerns positive association.

LEMMA 3.5. *The system $(\sigma_v, v \in \mathbb{Z}^2)$, described in Section 2.1, is positively associated. That is, for all increasing (in terms of the $\sigma$ variables) events $A$ and $B$, $P^{(h)}(A \cap B) \geq P^{(h)}(A)P^{(h)}(B)$.*

PROOF. The random variables $(X_i, i \in I)$ are independent $\{0, 1, \ldots, k\}$-valued random variables and hence, by FKG (or, in this special case, Harris' inequality for positive association), positively associated. Since the $\sigma$ variables are increasing functions of the $X$ variables, the statement of the lemma follows. □

REMARK. Note that the $\sigma_v, v \in \mathbb{Z}^2$ do not necessarily satisfy the strong FKG condition.

3.4. *RSW properties.* As stated in the Introduction, Bollobás and Riordan (see Theorem 4.1 in [8]) obtained a new RSW-like result for the Voronoi percolation model. The conclusion of their RSW theorem is weaker than that of the classical RSW theorem, but its proof is more robust: it does not (like the proof of "classical" RSW) involve conditioning on the lowest crossing. It works, as they pointed out, not only for the Voronoi model, but also for a large class of percolation models. In fact, the conditions are as follows (see [9] and Section 4.3 in [4]):

(a) crossings of rectangles must be defined in terms of "geometric paths" in such a way that (e.g.) horizontal and vertical crossings meet (this enables the often-used tool of pasting together paths to be used);

(b) certain increasing events (in particular, events of the form that there is a + path between two given sets of vertices) must be positively correlated;

(c) the distribution of the random field $(\sigma_v, v \in \mathbb{Z}^2)$ should be invariant under the symmetries of $\mathbb{Z}^2$;

(d) finally, certain mixing properties are needed.

The model in Theorem 2.2 satisfies the above conditions: as for (a), these are simply well-known properties for percolation on the square lattice and its matching lattice, and have nothing to do with the distribution $P^{(h)}$. As for (b) and (c), these are taken care of by Lemma 3.5 and by property (iv) in Section 2.1, respectively. Finally, as for (d), the following property (here formulated in our notation) is more than enough (see Remark 4.5 in [4]): for each $\varepsilon > 0$, there is an $l$ such that for all $k > l$, all $k$ by $2k$ rectangles $R_1$ and $R_2$ that have distance larger than $k/100$ to each other and all events $A$ and $B$ that are defined in terms of the random variables $(\sigma_v, v \in R_1)$ and the random variables $(\sigma_v, v \in R_2)$, respectively, $|P^{(h)}(A \cap B) - P^{(h)}(A)P^{(h)}(B)| < \varepsilon$. For our model, this is immediately guaranteed by Lemma 3.4. Hence, our model belongs to the class of models mentioned above.

For this class of models, the Bollobás–Riordan RSW-like theorem says that if the lim inf, as $s \to \infty$, of the probability of having a horizontal crossing of the box $[0, s] \times [0, s]$ is positive, then for every $\rho > 0$, the lim $\sup_{s \to \infty}$ of the probability of a horizontal crossing of the box $[0, \rho s] \times [0, s]$ is positive.

It is pointed out in [4] that small modifications of the proof of Theorem 4.1 in [8] in fact give the stronger result (for the same class of models as described above) that if for *some* $\rho > 0$, the lim $\sup_{s \to \infty}$ of the probability that there is a horizontal crossing of the box $[0, \rho s] \times [0, s]$ is positive, then this holds for *all* $\rho > 0$. (Note the occurrence of lim sup and lim inf.) Or, equivalently, if for some $\rho$, $\lim_{s \to \infty}$ of the probability that there is a horizontal crossing of the box $[0, \rho s] \times [0, s]$ equals 0, then this limit equals 0 for every $\rho > 0$. As remarked above, our current percolation model satisfies the required properties. So we get the aforementioned RSW result. Before we state this explicitly, we introduce the following notation. Let $H(n, m)$ [resp. $V(n, m)$] denote the event that there is a horizontal (resp. vertical) + crossing of the box $[0, n] \times [0, m]$. Further, let $H^{-*}(n, m)$ and $V^{-*}(n, m)$ be the analogs of $H(n, m)$ and $V(n, m)$ for $-$ crossings in the $*$ lattice. In this notation, the above mentioned RSW-like statement is as follows.

LEMMA 3.6. (a) *If*

$$\lim_{n \to \infty} P^{(h)}(H(\rho n, n)) = 0 \qquad \text{for some } \rho > 0,$$

*then*

$$\lim_{n \to \infty} P^{(h)}(H(\rho n, n)) = 0 \qquad \text{for all } \rho > 0.$$

(b) *The analogous result, with $H$ replaced by $H^{-*}$, also holds.*

Note that since a box either has a horizontal $+$ crossing or a vertical $-*$ crossing (and using rotation symmetry), we have that for each $k$ and $l$, $P^{(h)}(H(k,l)) = 1 - P^{(h)}(H^{-*}(l,k))$. Combining this with Lemma 3.6 immediately gives the following corollary.

COROLLARY 3.7. (a) *If*
$$\lim_{n \to \infty} P^{(h)}(H(\rho n, n)) = 1 \qquad \text{for some } \rho > 0,$$
*then*
$$\lim_{n \to \infty} P^{(h)}(H(\rho n, n)) = 1 \qquad \text{for all } \rho > 0.$$

(b) *The analogous result, with $H$ replaced by $H^{-*}$, also holds.*

3.5. *Finite-size criterion.*

LEMMA 3.8. *There is an $\hat{\varepsilon} > 0$ and an integer $\hat{N}$ such that for all $N \geq \hat{N}$, the following hold:*

(19) \qquad (a) \quad *if $P^{(h)}(V(3N, N)) < \hat{\varepsilon}$,*

*then the distribution of $|C^+|$ has an exponential tail;*

(20) \qquad (b) \quad *if $P^{(h)}(V^{-*}(3N, N)) < \hat{\varepsilon}$,*

*then the distribution of $|C^{-*}|$ has an exponential tail.*

PROOF. The proof below follows the main line of reasoning in the proof of the analogous well-known result for ordinary percolation (see [19]). Let $N$ and $\varepsilon$ be such that $P^{(h)}(V(3N, N)) < \varepsilon$. Cover $\mathbb{Z}^2$ by squares
$$Q_N(x) := Nx + [0, N]^2, \qquad x \in \mathbb{Z}^2.$$
We will often simply write $Q_N$ for $Q_N(\mathbf{0})$.

We say that an $x \in \mathbb{Z}^2$ is *good* if $Q_N(x)$ contains a vertex of $C^+$. A set $W \subset \mathbb{Z}^2$ is called good if every $x \in W$ is good. Let $S$ denote the set of good vertices. From the definition of "good," it is easy to see that $S$ is a connected subset of the square lattice and that $0 \in S$ unless $C^+ = \varnothing$ (in which case, also $S = \varnothing$). It is also clear that $|S| \geq |C^+|/|Q_N|$ and hence that

(21) \qquad $P^{(h)}(|C^+| \geq n) \leq P^{(h)}\left(|S| \geq \dfrac{n}{|Q_N|}\right), \qquad n = 1, 2, \ldots.$

Let, for $x \in \mathbb{Z}^2$, $R_1(x)$ denote the $3N \times N$ rectangle "north" of $Q_N(x)$. More precisely,
$$R_1(x) := Nx + [-N, 2N] \times [N, 2N].$$

Similarly, let $R_2(x)$ be the $3N \times N$ rectangle south of $Q_N(x)$ and let $R_3(x)$ and $R_4(x)$ be the $N \times 3N$ rectangles east, respectively west, of $Q_N(x)$.

Define, for each $x \in \mathbb{Z}^2$, the following event (where "easy" stands for "vertical" in the case of a $3N \times N$ rectangle and for "horizontal" in the case of an $N \times 3N$ rectangle):

$$A_x := \{\exists i \in \{1, \ldots, 4\} \text{ such that } R_i(x) \text{ has a } + \text{ crossing in the easy direction}\}.$$

It is standard (and easy to check) that for all (except a finite number, say $C_1$) $x \in \mathbb{Z}^2$, the following inclusion of events holds:

(22) $$\{x \text{ is good}\} \subset A_x.$$

Let $R(x) = \bigcup_{i=1}^{4} R_i(x)$. Recall the definition of "$l$ determined" in Section 3.2. We trivially have

(23) $$A_x \subset B_x,$$

where

$$B_x := (A_x \cap \{R(x) \text{ is } N \text{ determined}\}) \cup \{R(x) \text{ is not } N \text{ determined}\}.$$

We then get

(24)
$$P^{(h)}(B_x) \leq P^{(h)}(A_x) + P^{(h)}[R(x) \text{ is not } N \text{ determined}]$$
$$\leq 4\varepsilon + |R(0)| \max_x P^{(h)}(x \text{ is not } N \text{ determined})$$
$$\leq 4\varepsilon + C_2 N^2 \frac{C_0}{N^{2+\gamma}}$$
$$\leq 4\varepsilon + C_3(N), \qquad \text{where } C_3(N) \to 0 \text{ as } N \to \infty$$

and where the first inequality is trivial, the second follows from our choice of $N$ and $\varepsilon$, the third follows from (14) and $C_2$ is a constant.

Let $\alpha$ be as in property (iii) in Section 2.1. It is easy to see that there is a constant $C_4 = C_4(\alpha)$ such that for every finite set of vertices $x(1), \ldots, x(m)$ satisfying $\min_{1 \leq i < j \leq m} \|x(i) - x(j)\| > C_4(\alpha)$, the events $B_{x(i)}$, $1 \leq i \leq m$, are independent.

From this [and (22)–(24)], it follows easily that there exist a $C_5(\alpha)$ and $C_6(\alpha)$ such that for every finite set of vertices $W$,

(25)
$$P^{(h)}(W \text{ is good }) \leq (4\varepsilon + C_3(N))^{\lfloor(|W|-C_1)/C_5(\alpha)\rfloor}$$
$$\leq^* (4\varepsilon + C_3(N))^{|W|/C_6(\alpha)},$$

where the mark * in the last inequality means that inequality holds for all values of $|W|$ that are sufficiently large.

We now apply this to (21). To do this, note that if $|S| \geq n/|Q_N|$, then there is a good lattice animal $W$ of size $\lfloor \frac{n}{|Q_N|} \rfloor$. (A lattice animal is a connected set of vertices containing 0.) Using this, (21), (25) and the fact that there is a constant $C_7$ such that the number of lattice sites of size $k$ is at most $C_7^k$, we get

$$\begin{aligned}(26)\quad P^{(h)}(|C^+| \geq n) &\leq C_7^{\lfloor n/|Q_N|\rfloor}(4\varepsilon + C_3(N))^{\lfloor n/|Q_N|\rfloor/C_6(\alpha)} \\ &\leq C_8(\varepsilon, N)[(C_7(4\varepsilon + C_3(N))^{1/C_6(\alpha)})^{1/|Q_N|}]^n.\end{aligned}$$

Now, take $\hat{\varepsilon}$ and $\hat{N}$ such that $C_7(4\hat{\varepsilon}+C_3(\hat{N}))^{1/C_6(\alpha)} < 1$ for all $N \geq \hat{N}$ [which can be done since $C_3(N) \to 0$ as $N \to \infty$]. From (26), it follows that for this choice of $\hat{\varepsilon}$ and $\hat{N}$, the statement in part (a) of Lemma 3.8 holds. By exactly the same arguments (and, if necessary, by decreasing, resp. increasing, the values of $\hat{\varepsilon}$ and $\hat{N}$ obtained above), part (b) also follows. $\square$

From the above lemma, we easily get the following.

COROLLARY 3.9. *Let $\hat{\varepsilon}$ and $\hat{N}$ be as in Lemma 3.8.*

(a) *If there is an $n \geq \hat{N}$ with $P^{(h)}(V(3n,n)) < \hat{\varepsilon}$, then $P^{(h)}(|C^{-*}| = \infty) > 0$.*

(b) *If there is an $n \geq \hat{N}$ with $P^{(h)}(V^{-*}(3n,n)) < \hat{\varepsilon}$, then $P^{(h)}(|C^+| = \infty) > 0$.*

PROOF. We only prove part (a) here; the proof of (b) is completely analogous. If the condition of Corollary 3.9 holds, then by Lemma 3.8, the distribution of $|C^+|$ has an exponential tail. Exactly as in the Peierls argument in ordinary percolation (see, e.g., [13]), this implies that the probability that there is a + circuit having 0 in its interior is less than 1 and hence that $P^{(h)}(|C^{-*}| = \infty) > 0$. $\square$

LEMMA 3.10. *If $P^{(h)}(|C^+| = \infty) > 0$, then $P^{(h)}(|C^{-*}| = \infty) = 0$.*

PROOF. There are various standard ways to prove this. One is as follows. The law of $(\sigma_v, v \in \mathbb{Z}^2)$, is positively associated (by Lemma 3.5), translation invariant, invariant under horizontal axis reflection and vertical axis reflection [property (iv) in Section 2.1] and mixing (in the ergodic-theoretic sense, w.r.t. horizontal translations as well as to vertical translations). The last follows from Lemma 3.4. Hence, by the main result in [11], $P^{(h)}(|C^{-*}| = \infty) = 0$. $\square$

**4. Proof of Theorem 2.2.** We use the notation $\theta(h)$ for $P^{(h)}(|C^+| = \infty)$ and $\theta^{-*}(h)$ for $P^{(h)}(|C^{-*}| = \infty)$. Let

$$h_c := \sup\{h : \theta(h) = 0\}.$$

It is quite easy to see that $h_c < \infty$. Take $n \geq \hat{N}$ with $\hat{N}$ defined as in Lemma 3.8. From properties (b), (i) and (iii) in Section 2.1, it follows that for all $v \in \mathbb{Z}^2$, $P^{(h)}(\sigma_v = +1) \to 1$ as $h \to \infty$ and hence that $P^{(h)}(H(3n,n)) \to 1$ as $h \to \infty$, which is equivalent to $P^{(h)}(V^{-*}(3n,n)) \to 0$ as $h \to \infty$. So, there is an $h$ such that $P^{(h)}(V^{-*}(3n,n)) < \hat{\varepsilon}$ with $\hat{\varepsilon}$ as in Corollary 3.9. By part (b) of that corollary, $\theta(h) > 0$ for such $h$. Hence, we indeed have that $h_c < \infty$. Using analogous arguments, it follows that $h_c > -\infty$.

PROOF OF THEOREM 2.2. We now start with the proof of part (a) of Theorem 2.2, where we will use the following notation. $B(n)$ denotes the square $[-n,n]^2$ and $\partial B(n)$ its boundary [the set of all vertices $v$ that are not in $B(n)$, but for which there is a $w \in B(n)$ with $\|v - w\| = 1$]. For $n \leq m$, $A(n,m)$ denotes the annulus $B(m) \setminus B(n)$. For $v \in \mathbb{Z}^2$ and $n \in N$, $B(v;n)$ will denote the set $B(n)$ shifted by $v$.

Let $h$ be larger than the above-defined $h_c$. So, $P^{(h)}(|C^+| = \infty) > 0$. We will first show that

$$P^{(h)}(H(n,n)) \to 1 \quad \text{as } n \to \infty. \tag{27}$$

This is done in a quite standard way. Let $\delta > 0$. Take $K$ sufficiently large that

$$P^{(h)}(B(K) \leftrightarrow \infty) > 1 - \delta. \tag{28}$$

By Lemma 3.10, we can take $N > K$ so large that

$$P^{(h)}(\exists \text{ a } + \text{ circuit in } A(K,N) \text{ surrounding } B(K)) > 1 - \delta. \tag{29}$$

For all $n \geq N$, the following holds. First, by (28), we have, of course, that $P^{(h)}(B(K) \leftrightarrow \partial B(n)) > 1 - \delta$. Since our model has the positive association property (see, Lemma 3.5), we can apply the usual "square root trick" (see, e.g., [13]), which gives that $P^{(h)}(B(K) \leftrightarrow r(B(n))) > 1 - \delta^{1/4}$, where $r(B(n))$ denotes the right-hand side $\{n\} \times [-n,n]$ of $B(n)$. By this and its analog for the left side $l(B(n))$ of $B(n)$, together with (29) (and again positive association), we get, for all $n \geq N$,

$$\begin{aligned} &P^{(h)}(H(n,n)) \\ &\geq P^{(h)}(B(K) \leftrightarrow r(B(n)), B(K) \leftrightarrow l(B(n)), + \text{ circuit in } A(K,n)) \\ &\geq (1 - \delta^{1/4})^2(1 - \delta). \end{aligned} \tag{30}$$

Since we can take $\delta$ arbitrary small, (27) follows.

Application of Corollary 3.7 now gives that $P^{(h)}(H(3n,n)) \to 1$ as $n \to \infty$ and hence that $P^{(h)}(V^{-*}(3n,n)) \to 0$ as $n \to \infty$.

Finally, by part (b) of Lemma 3.8, this implies that the distribution of $|C^{-*}|$ has an exponential tail. This completes the proof of Theorem 2.2(a).

It is important to note that part (b) of the theorem cannot simply be concluded by replacing "+" by "−∗" (and vice versa) in the arguments above. The problem is that our definition of $h_c$ in the beginning of the proof is "asymmetric." If we could show that the above defined $h_c$ is equal to $\inf\{h : P^{(h)}(|C^{-*}| = \infty) = 0\}$ or, equivalently [since we already know, by Lemma 3.10, that there is no $h$ for which both $\theta(h) > 0$ and $\theta^{-*}(h) > 0$], that $\theta^{-*}(h) > 0$ for all $h < h_c$, we would be able to conclude (b) by exchanging $+$ and $-*$ in the arguments of (a). Below, it will be shown, using the approximate zero-one laws in Section 3.1, that, indeed, $\theta^{-*}(h) > 0$ for all $h < h_c$.

*Proof of* (b). Suppose there is an $h_1 < h_c$ with $\theta^{-*}(h_1) = 0$. We will show that this leads to a contradiction. Let $h_2 \in (h_1, h_c)$. Then, for all $h \in [h_1, h_2]$, by monotonicity [see properties (a) and (i) in Section 2.1], $\theta(h) = \theta^{-*}(h) = 0$. Let $H(n, m)$ and $H^{-*}(n, m)$ be the box-crossing events defined in Section 3.4. Since $\theta^+ \equiv 0$ on $[h_1, h_2]$, we have, by Corollary 3.9 (b) [noting that $P^{(h)}(V^{-*}(3n, n)) = 1 - P^{(h)}(H(3n, n))$], that

$$(31) \qquad \forall h \in [h_1, h_2] \ \forall n \geq \hat{N} \qquad P^{(h)}(H(3n, n)) < 1 - \hat{\varepsilon}$$

with $\hat{\varepsilon}$ and $\hat{N}$ as in Lemma 3.8.

On the other hand, $P^{(h_1)}(H(n, 3n)) = P^{(h_1)}(V(3n, n))$, which [again by Corollary 3.9 and because $\theta^{-*}(h_1) = 0$] is at least $\hat{\varepsilon}$ for all $n \geq \hat{N}$. Hence, by Lemma 3.6,

$$\limsup_{n \to \infty} P^{(h_1)}(H(3n, n)) > 0.$$

Using this, monotonicity and (31), it follows straightforwardly that there is a $\delta \in (0, 1)$ and an infinite sequence $n_1 < n_2 < n_3 < \cdots$ such that

$$(32) \qquad P^{(h)}(H(3n_i, n_i)) \in (\delta, 1 - \delta) \qquad \text{for all } i \text{ and all } h \in [h_1, h_2].$$

To reach a contradiction, we will show that the sequence $(\varepsilon_n)$, defined by

$$\varepsilon_n := \sup_{j \in I, h \in [h_1, h_2]} P^{(h)}((H(3n, n))_j),$$

satisfies

$$(33) \qquad \varepsilon_n \to 0 \qquad \text{as } n \to \infty,$$

where, as before (Section 3.1), $A_j$ denotes the event that $j$ is an internal pivotal index for the event $A$.

REMARK. It is important to note that, here, we do not (as in ordinary percolation and in Higuchi's treatment) consider pivotality in terms of the vertices of the lattice (the indices of the $\sigma$ variables), but in terms of the indices of the underlying $X$ variables [i.e., in the special case of the Ising model, the space-time variables $Y_v(t)$ in Section 2.2.1, which control the updates in the dynamics].

We will first show that (32) and (33) indeed give a contradiction. By Theorem 3.3 we have, for all $i = 1, 2, \ldots$,

$$P^{(h_1)}(H(3n_i, n_i))(1 - P^{(h_2)}(H(3n_i, n_i))) \leq (\varepsilon_{n_i})^{c(h_1, h_2)(h_2 - h_1)/K_2}$$

with $\varepsilon_{n_i}$ as defined above.

By (33), the right-hand side in this last inequality goes to 0 as $n \to \infty$. However, for all $i$, the left-hand side is at least $\delta^2$ by (32)—a contradiction.

So, part (b) of the theorem is proved once we prove (33), which we will do now. In the following, $X$ stands for the collection of random variables $(X_i, i \in I)$.

Note that, by the definition of internal pivotal,

$$(34) \qquad P^{(h)}((H(3n, n))_j) = P(X \in H(3n, n), X^{(j)} \notin H(3n, n)),$$

where $X^{(j)}$ is the element of $\{0, \ldots, k\}^I$ that satisfies $X_i^{(j)} = X_i$ for all $i \neq j$ and $X_j^{(j)} = 0$.

Now, recall that for each $v \in \mathbb{Z}^2$, we have the sequence $i_1(v), i_2(v), \ldots$ introduced in property (ii) of Section 2.1. We will use the following terminology. If $i_m(v) = j$, we say that $j$ has rank $m$ for $v$. If $j$ does not occur at all in the sequence $i_1(v), i_2(v), \ldots$, we say that the rank of $j$ for $v$ is infinite. The rank of $j$ for $v$ will be denoted by $r_v(j)$. Suppose that $r_v(j) = m$. Then we say that $v$ *needs* $j$ if $(X_{i_1(v)}, X_{i_2(v)}, \ldots, X_{i_{m-1}(v)})$ does not determine $\sigma_v$.

Let $v$ be a vertex in the box $[0, 3n] \times [0, n]$. We use the notation $H(3n, n; v)$ for the event that $v$ is on a horizontal + crossing of that box. Using the terminology and observation above, we have that the right-hand side of (34) is at most

$$P(\exists v \in \mathbb{Z}^2 \text{ such that } X \in H(3n, n; v), \text{ but } X^{(j)} \notin H(3n, n; v))$$

$$\leq \sum_{v \in \mathbb{Z}^2} P(X \in H(3n, n; v), X^{(j)} \notin H(3n, n; v))$$

(35)

$$\leq \sum_{v \in \mathbb{Z}^2} P^{(h)}(H(3n, n; v), v \text{ needs } j)$$

$$\leq \sum_{v \in \mathbb{Z}^2} \min(P^{(h)}(H(3n, n; v)), P^{(h)}(v \text{ needs } j)).$$

Further, note that if $v$ is on a horizontal + crossing of the rectangle $[0, 3n] \times [0, n]$, there must be a + path from $v$ to $\partial B(v; n)$. By this, and translation invariance [property (iv) of Section 2.1], the first of the two probabilities in the expression in the summand in the last line of (35) [i.e., $P^{(h)}(H(3n, n; v))$] is at most $P^{(h)}(0 \leftrightarrow \partial B(n))$, which, by monotonicity, is, of course, at most $P^{(h_2)}(0 \leftrightarrow \partial B(n))$. Let us denote this last probability by $f(n)$. Also, note that property (ii) of Section 2.1 states that

$$P^{(h)}(v \text{ needs } j) \leq \frac{C_0}{(r_v(j)-1)^{2+\gamma}}.$$

These considerations imply that the last line of (35) is, for each positive integer $K$, at most

$$(36) \quad f(n) \times |\{v \in \mathbb{Z}^2 : r_v(j) \leq K\}| + C_0 \sum_{k=K}^{\infty} \frac{|\{v : r_v(j) = k\}|}{(k-1)^{2+\gamma}}.$$

Consider the set in the first term in (36). Let $u$ and $w$ be two vertices which both belong to this set. That is, $r_u(j) \leq K$ and $r_w(j) \leq K$ hold and hence the sets $\{i_1(u), \ldots, i_K(u)\}$ and $\{i_1(w), \ldots, i_K(w)\}$ have nonempty intersection. It follows from property (iii) of Section 2.1 that $\|v - w\|$ is at most $K/\alpha$. Hence, the set under consideration has diameter $\leq K/\alpha$. The cardinality of this set therefore satisfies

$$(37) \quad |\{v \in \mathbb{Z}^2 : r_v(j) \leq K\}| \leq \frac{C_9 K^2}{\alpha^2}$$

for some constant $C_9$.

From this [and using the fact that $(k-1)^{2+\gamma}$ is decreasing in $k$], it is easy to see that the sum in (36) satisfies

$$(38) \quad \sum_{k=K}^{\infty} \frac{|\{v : r_v(j) = k\}|}{(k-1)^{2+\gamma}} \leq \frac{C_{10}}{\alpha^2 K^\gamma} + \sum_{k=K+1}^{\infty} \frac{C_{10}}{\alpha^2 k^{1+\gamma}}$$

for some constant $C_{10}$.

Note that in (36), we are free to choose $K$. In the following, we let $K(n)$ be the largest integer $k$ for which

$$\frac{C_9 k^2}{\alpha^2} \leq \frac{1}{\sqrt{f(n)}}.$$

Taking together (34)–(38) we get, choosing $K = K(n)$ in (36),

$$(39) \quad P^{(h)}((H(3n, n))_j) \leq f(n) \frac{1}{\sqrt{f(n)}} + \frac{C_{10}}{\alpha^2 K(n)^\gamma} + \sum_{k=K(n)+1}^{\infty} \frac{C_{10}}{\alpha^2 k^{1+\gamma}}.$$

Note that the right-hand side of (39) does not depend on $j$ and $h$, and [since $f(n) \to 0$ as $n \to \infty$, $\gamma > 0$ and $K(n) \to \infty$ as $n \to \infty$] goes to 0 as $n \to \infty$. This proves (33) and thus completes the proof of the first statement in part (b) of the theorem. The second statement of part (b) now follows in exactly the same way as its analog in (a). □

**Acknowledgments.** My interest in some of the problems in this paper was raised by communication with Frank Redig and Federico Camia. I also thank Raphaël Rossignol for valuable discussions concerning Section 3, and the referee for many detailed comments.

## REFERENCES

[1] AIZENMAN, M. and BARSKY, D. (1987). Sharpness of the phase transition in percolation models. *Comm. Math. Phys.* **108** 489–526. MR0874906

[2] AIZENMAN, M., BARSKY, D. and FERNANDEZ, R. (1987). The phase transition in a general class of Ising-type models is sharp. *J. Statist. Phys.* **47** 343–374. MR0894398

[3] BALINT, A., CAMIA, F. and MEESTER, R. (2007). Sharp phase transition and critical behaviour in $2D$ divide and colour models. *Stoch. Proc. Appl.* To appear.

[4] VAN DEN BERG, J., BROUWER, R. and VÁGVÖLGYI, B. (2008). Box-crossings and continuity results for self-destructive percolation in the plane. *In and Out of Equilibrium II* (V. Sidoravicius and M.-E. Vares, eds.). *Progr. Probab.* To appear.

[5] VAN DEN BERG, J. and KESTEN, H. (1985). Inequalities with applications to percolation and reliability. *J. Appl. Probab.* **22** 556–569. MR0799280

[6] VAN DEN BERG, J. and STEIF, J. E. (1999). On the existence and nonexistence of finitary codings for a class of random fields. *Ann. Probab.* **27** 1501–1522. MR1733157

[7] BOLLOBÁS, B. and RIORDAN, O. (2006). A short proof of the Harris–Kesten theorem. *Bull. London Math. Soc.* **38** 470–484. MR2239042

[8] BOLLOBÁS, B. and RIORDAN, O. (2006). The critical probability for random Voronoi percolation in the plane is $1/2$. *Probab. Theory Related Fields* **136** 417–468. MR2257131

[9] BOLLOBÁS, B. and RIORDAN, O. (2006). Sharp thresholds and percolation in the plane. *Random Structures Algorithms* **29** 524–548. MR2268234

[10] FRIEDGUT, E. and KALAI, G. (1996). Every monotone graph property has a sharp threshold. *Proc. Amer. Math. Soc.* **124** 2993–3002. MR1371123

[11] GANDOLFI, A., KEANE, M. and RUSSO, L. (1988). On the uniqueness of the infinite occupied cluster in dependent two-dimensional site percolation. *Ann. Probab.* **16** 1147–1157. MR0942759

[12] GRAHAM, B. T. and GRIMMETT, G. R. (2006). Influence and sharp-threshold theorems for monotonic measures. *Ann. Probab.* **34** 1726–1745. MR2271479

[13] GRIMMETT, G. R. (1999). *Percolation*, 2nd ed. Springer, Berlin. MR1707339

[14] HARRIS, T. E. (1960). A lower bound for the critical probability in a certain percolation process. *Proc. Cambridge Philos. Soc.* **56** 13–20. MR0115221

[15] HIGUCHI, Y. (1993). Coexistence of infinite (*)-clusters II. Ising percolation in two dimensions. *Probab. Theory Related Fields* **97** 1–33. MR1240714

[16] HIGUCHI, Y. (1993). A sharp transition for the two-dimensional Ising percolation. *Probab. Theory Related Fields* **97** 489–514. MR1246977

[17] KAHN, J., KALAI, G. and LINIAL, N. (1988). The influence of variables on Boolean functions. In *Proc. 29-th Annual Symposium on Foundations of Computer Science* 68–80. Computer Society Press.

[18] KESTEN, H. (1980). The critical probability of bond percolation on the square lattice equals $\frac{1}{2}$. *Comm. Math. Phys.* **74** 41–59. MR0575895

[19] KESTEN, H. (1981). Analyticity properties and power law estimates of functions in percolation theory. *J. Statist. Phys.* **25** 717–756. MR0633715

[20] KESTEN, H. (1987). Scaling relations for 2D percolation. *Comm. Math. Phys.* **109** 109–156. MR0879034

[21] LEBOWITZ, J. L. (1974). GHS and other inequalities. *Comm. Math. Phys.* **35** 87–92. MR0339738

[22] MARTINELLI, F. and OLIVIERI, E. (1994). Approach to equilibrium of Glauber dynamics in the one phase region. I. The attractive case. *Comm. Math. Phys.* **61** 447–486. MR1269387
[23] MENSHIKOV, M. V. (1986). Coincidence of critical points in percolation problems. *Soviet Math. Dokl.* **33** 856–859. MR0852458
[24] PROPP, J. G. and WILSON, D. B. (1996). Exact sampling with coupled Markov chains and applications to Statistical Mechanics. *Random Structures Algorithms* **9** 223–252. MR1611693
[25] ROSSIGNOL, R. (2006). Threshold for monotone symmetric properties through a logarithmic Sobolev inequality. *Ann. Probab.* **34** 1707–1725. MR2271478
[26] ROSSIGNOL, R. (2008). Threshold phenomena on product spaces: BKKKL revisited (once more). *Electron. Comm. Probab.* **13** 35–44. MR2372835
[27] RUSSO, L. (1978). A note on percolation. *Z. Wahrsch. Verw. Gebiete* **43** 39–48. MR0488383
[28] RUSSO, L. (1981). On the critical percolation probabilities. *Z. Wahrsch. Verw. Gebiete* **56** 229–237. MR0618273
[29] RUSSO, L. (1982). An approximate zero-one law. *Z. Wahrsch. Verw. Gebiete* **61** 129–139. MR0671248
[30] SEYMOUR, P. D. and WELSH, D. J. A. (1978). Percolation probabilities on the square lattice. In *Advances in Graph Theory* (B. Bollobás ed.) 227–245. *Ann. Discrete Math.* **3**. North-Holland, Amsterdam. MR0494572
[31] TALAGRAND, M. (1994). On Russo's approximate zero-one law. *Ann. Probab.* **22** 1576–1587. MR1303654
[32] WIERMAN, J. C. (1981). Bond percolation on honeycomb and triangular lattices. *Adv. in Appl. Probab.* **13** 293–313. MR0612205

CWI
KRUISLAAN 413
1098 SJ AMSTERDAM
THE NETHERLANDS
E-MAIL: J.van.den.Berg@cwi.nl